\documentclass[reqno]{amsart}
\newtheorem{thm}{Theorem}

\theoremstyle{definition}

\newtheorem{ex}[thm]{Example}

\providecommand{\abs}[1]{\lvert#1\rvert}

\providecommand{\norm}[1]{\lVert#1\rVert}
\begin{document}

\title[Bayesian predictive inference]{Some cautionary tales about\\Bayesian predictive inference}

\author{Emanuela Dreassi}
\address{Emanuela Dreassi, Dipartimento di Statistica, Informatica, Applicazioni ``G. Parenti'', Universit\`a di Firenze, viale Morgagni 59, 50134 Firenze, Italy}
\email{emanuela.dreassi@unifi.it}
\author{Fabrizio Leisen}
\address{Fabrizio Leisen, Department of Mathematics, King’s College, Strand WC2R 2LS, London, UK} \email{fabrizio.leisen@gmail.com}
\author{Luca Pratelli}
\address{Luca Pratelli, Accademia Navale, viale Italia 72, 57100 Livorno,
Italy} \email{luca\_pratelli@marina.difesa.it}
\author{Pietro Rigo}
\address{Pietro Rigo (corresponding author), Dipartimento di Scienze Statistiche ``P. Fortunati'', Universit\`a di Bologna, via delle Belle Arti 41, 40126 Bologna, Italy}
\email{pietro.rigo@unibo.it}



\begin{abstract}
Two misunderstandings, frequently arising in Bayesian predictive inference, are discussed. The first deals with the data generating mechanism, while the second consists in overestimating the role played by asymptotic exchangeability. Some consequences of such misunderstandings are highlighted through examples.
\end{abstract}

\maketitle

\section{Introduction}\label{intro}

The only purpose of this note is to highlight, and possibly clarify, two common misunderstandings which frequently arise in Bayesian predictive inference. This goal is achieved through some discussion and a few related examples. These misunderstandings are not independent, but, to make the exposition clearer, we treat them separately.

We begin with some notation. Throughout,
$$X=(X_1,X_2,\ldots)$$
is a sequence of real random variables, meant as the data sequence. We denote by $\mathcal{B}$ the Borel $\sigma$-field on $\mathbb{R}$ and we assume that all random elements are defined on a common probability space, say $(\Omega,\mathcal{A},P)$. The {\em predictive distributions} of $X$ are
$$\alpha_0(\cdot)=P(X_1\in\cdot)\quad\text{and}\quad\alpha_n(\cdot)=P(X_{n+1}\in\cdot\mid X_1,\ldots,X_n).$$
Thus, $\alpha_0$ is the marginal distribution of $X_1$ while $\alpha_n$ is the conditional distribution of $X_{n+1}$ given $(X_1,\ldots,X_n)$.

We recall that a {\em random probability measure} (r.p.m.) is a function $\alpha$ on $\Omega\times\mathcal{B}$ such that $\alpha(\omega,\cdot)$ is a probability measure on $\mathcal{B}$, for fixed $\omega\in\Omega$, and $\alpha(\cdot,A)$ is a real random variable for fixed $A\in\mathcal{B}$. Obviously, a predictive distribution $\alpha_n$ is a special type of r.p.m. Another example of r.p.m. is the {\em empirical measure}
$$\mu_n=\frac{1}{n}\,\sum_{i=1}^n\delta_{X_i}$$
where $\delta_x$ denotes the unit mass at the point $x$.

By the Ionescu-Tulcea theorem (see e.g. \cite{STS25}) the probability distribution of the sequence $X$ is uniquely determined by its predictives $\alpha_n$. More importantly, apart from some mild measurability conditions, the $\alpha_n$ can be selected arbitrarily. Roughly speaking, to assign the probability distribution of $X$, it suffices to select (arbitrarily) the probability distribution of $X_1$, the conditional distribution of $X_2$ given $X_1$, the conditional distribution of $X_3$ given $(X_1,X_2)$, and so on.

In fact, the assignment of a sequence $(\alpha_n)$ of predictives is the fundamental step of a number of statistical procedures. The obvious example is Bayesian predictive inference \cite{STS25,FP}, which is also the reference of this note. In passing, however, we recall that the choice of $(\alpha_n)$ is basic in other frameworks, including species sampling models \cite{PIT}, predictive resampling \cite{FHW,FPJB}, machine learning \cite{DUT}, causal inference \cite{LDM}, and Dawid's prequential approach \cite{DV}.

\subsection{First misunderstanding}\label{xd56n}
There are two conflicting views regarding the data sequence $X$. Both are admissible but mixing the two can cause serious issues. The first misunderstanding, hereafter M1, arises when it is unclear, or at least vague, which point of view about $X$ is adopted.

According to the {\em first view}, the underlying probability measure $P$ is induced by the predictives $\alpha_n$ via Ionescu-Tulcea theorem. Precisely, the $X_n$ are the coordinate random variables on $(\Omega,\mathcal{A})=(\mathbb{R}^\infty,\mathcal{B}^\infty)$ and $P$ is the probability measure on $\mathcal{B}^\infty$ induced by $(\alpha_n)$. In other terms, once $(\alpha_n)$ has been selected, $P$ is uniquely determined and the $\alpha_n$ are the predictives of $X$ under $P$.

According to the {\em second view}, $P$ is unknown but such that $X$ satisfies some dependence structure, such as exchangeable, partially exchangeable, stationary, Markov, and so on. Moreover, the inferrer also selects a sequence $(\beta_n)$ of r.p.m.'s with $\beta_n$ measurable with respect to $\sigma(X_1,\ldots,X_n)$. The $\beta_n$ are not necessarily the predictives of $X$ under $P$. Nevertheless, they are often called ``predictive distributions" and they are used as if they were. More correctly, the $\beta_n$ should be regarded as estimators of some unknown (and possibly random) probability measure $\mu$. For instance, $X$ is i.i.d., $\mu$ is the unknown common distribution, and $\beta_n$ is an estimator of $\mu$ based on $(X_1,\ldots,X_n)$. Or else, to put it more rigorously, $X$ is exchangeable and $\beta_n$ is an estimator of $\mu$ where $\mu$ is a r.p.m. such that
\begin{gather*}
P(X\in A)=E\bigl\{\mu^\infty(A)\bigr\}\quad\text{for all }A\in\mathcal{B}^\infty.
\end{gather*}
Two final remarks are in order. Firstly, the $\beta_n$ should not be called ``predictive distributions". Secondly, if $X$ is assumed to be exchangeable and $\beta_n$ is regarded as an estimator of $\mu$, it is reasonable to require that $\beta_n$ is invariant under permutations of $(X_1,\ldots,X_n)$.

To make a few examples, Theorem 1 of \cite{HMW} and Theorem 7 of \cite{FHW} are meant according to the second view. The same happens for most results regarding Bayesian consistency (where $X$ is usually regarded as i.i.d. with an unknown probability distribution). Instead, Theorem 5 of \cite{FHW} and the results in \cite{BHW,PIT} refer to the first view. Similarly, all the authors' results are obtained in the first view; see e.g. \cite{AOPCID,AOP13,STS25,ECP26} and references therein.

To be fair, it should be noted that we usually adopt the first perspective about $X$. But this is of no importance. Our current goal is not to support one view or the other, but to distinguish them.

\subsection{Second misunderstanding}\label{bhy67c}
Say that $X$ is {\em asymptotically exchangeable} if
\begin{gather*}
(X_{n+1},X_{n+2},\ldots)\,\overset{d}\longrightarrow\, (Z_1,Z_2,\ldots),\quad\text{as }n\rightarrow\infty,
\end{gather*}
where $Z=(Z_1,Z_2,\ldots)$ is an exchangeable sequence and $\overset{d}\longrightarrow$ means convergence in distribution. In this case, $Z$ is called the exchangeable limit sequence of $X$.

By Lemma (8.2) of Aldous \cite{A}, for $X$ to be asymptotically exchangeable, it suffices that $\alpha_n\rightarrow\mu$, weakly a.s., for some r.p.m. $\mu$. Here, weak a.s. convergence is meant as
$$\alpha_n(\omega,\cdot)\overset{weakly}\longrightarrow\mu(\omega,\cdot),\quad\text{as }n\rightarrow\infty,\text{ for }P\text{-almost all }\omega\in\Omega.$$
In addition, if $\alpha_n\rightarrow\mu$ weakly a.s., the exchangeable limit sequence $Z$ is driven by $\mu$, i.e.
$$P(Z\in A)=E\bigl\{\mu^\infty(A)\bigr\}\quad\text{for all }A\in\mathcal{B}^\infty.$$
Aldous' result can be slightly improved, in the sense that weak convergence a.s. can be weakened into weak convergence in probability; see Theorem 2.2 of \cite{ECP26}. Generally, however, asymptotic exchangeability does not imply weak convergence of the $\alpha_n$, neither in probability; see Example 4.3 of \cite{ECP26} and Example \ref{cs22} below.

The second misunderstanding, hereafter M2, is the strong overestimation of the role played by asymptotic exchangeability.

To be more precise, for any sequence $Y$, write $Y\sim X$ if
\begin{gather*}
(Y_{n+1},Y_{n+2},\ldots)\overset{d}\longrightarrow Z\quad\text{and}\quad(X_{n+1},X_{n+2},\ldots)\overset{d}\longrightarrow Z
\end{gather*}
for some exchangeable sequence $Z$. Moreover, suppose that $X$ is asymptotically exchangeable and $\mu$ is the r.p.m. such that $P(Z\in\cdot)=E\bigl\{\mu^\infty(\cdot)\bigr\}$, where $Z$ is the exchangeable limit sequence of $X$.

Then, asymptotic exchangeability is often understood as if, for large $n$, the shifted sequence $(X_{n+1},X_{n+2},\ldots)$ would satisfy some version of de Finetti's representation theorem. In other terms, $(X_{n+1},X_{n+2},\ldots)$ is often regarded as i.i.d., conditionally on $\mu$, with common distribution $\mu$. See e.g. \cite{BHW,CSRS,FHW,FPJB,FP}. Or else, sometimes, all sequences in the class $\{Y:Y\sim X\}$ are seen as equivalent for inferential purposes. Taking this interpretation one step further, if $X\sim Y$, then $X$ and $Y$ should yield the same asymptotic inference on the object of interest (whatever the latter is meant).

The previous interpretations are excessive, even admitting a certain degree of approximation. Asymptotic exchangeability is a weak condition, and the class $\bigl\{Y:Y\sim X\bigr\}$ is very large and heterogeneous. For instance, it may be that $X\sim Y$ and yet $X$ and $Y$ have singular probability distributions, in the sense that $P(X\in A)=0$ and $P(Y\in A)=1$ for some $A\in\mathcal{B}^\infty$; see Example \ref{cs22}. Note also that $Z\sim X$, where $Z$ is the exchangeable limit sequence of $X$. Hence, according to the previous interpretations, $X$ and $Z$ should be inferentially equivalent. This is rather stretched. As an example, suppose the object of inference is a function of the data sequence, say $\theta=f(X)$ where $f$ is a (known) measurable function on $\mathbb{R}^\infty$. To be concrete, suppose $\theta=\lim_n\frac{1}{n}\,\sum_{i=1}^nX_i$ (where the a.s. existence of the limit is an assumption). Then, $X$ and $Z$ may yield very different inferences, on $\theta=f(X)$ and $\theta=f(Z)$ respectively, even though $X\sim Z$; see Example \ref{bg67j}.

\section{Examples}

The consequences of M1 and M2 may be highlighted via some examples. In the heading of each example, we write M1 or M2 to indicate which misunderstanding the example is dealing with.

While obvious, the next example helps to make M1 clear.

\begin{ex}\textbf{(M1; Classical bootstrap)}\label{emp}
Having the bootstrap in mind, it is tempting to let $\alpha_n=\mu_n$ for all $n\ge 1$ where $\mu_n$ is the empirical measure. In the first view, $P$ is the probability measure on $\mathcal{B}^\infty$ induced by $(\alpha_0,\mu_1,\mu_2,\ldots)$. Hence, under $P$, the data sequence reduces to
$$X_1=X_2=X_3=\ldots$$
with $X_1$ distributed according to $\alpha_0$. In other terms, in the first view, $\alpha_n=\mu_n$ for all $n\ge 1$ is a very bad choice. On the contrary, this choice makes sense in the second view. For instance, suppose $X$ is assumed to be exchangeable and take $\beta_n=\mu_n$. Then, $\beta_n\rightarrow\mu$ weakly a.s. where $\mu$ is the r.p.m. such that $P(X\in\cdot)=E\bigl\{\mu^\infty(\cdot)\bigr\}$. Hence, $\beta_n$ is a reasonable and consistent estimator of $\mu$.
\end{ex}

In a sense, the next example is analogous to the previous one. But, at the same time, it is quite different.

\begin{ex}\textbf{(M1; Hill's predictive distributions)}\label{xf6yh8}
For each $n\ge 0$, let $\alpha_n$ be non-atomic and supported by $[0,1]$. Suppose also that
$$\alpha_n\Bigl[\bigl(X(j,n),\,X(j+1,n)\bigr)\Bigr]=\frac{1}{n+1}\quad\text{for }n\ge 1\text{ and }j=0,1,\ldots,n,$$
where $X(0,n)=0$, $X(n+1,n)=1$ and $(X(1,n),\ldots,X(n,n))$ are the order statistics based on $(X_1,\ldots,X_n)$. In some respects, $\alpha_n$ is analogous to $\mu_n$. Nevertheless,
$$\mu_n(A_n)=1\quad\text{while}\quad\alpha_n(A_n)=0\quad\text{where }A_n=\{X_1,\ldots,X_n\}.$$
The $\alpha_n$ have been introduced and motivated by Hill in \cite{HILL}. In \cite{BHW}, it is noted that the $\alpha_n$ are connected to conformal prediction. Assuming $\alpha_n$ uniform on each interval $\bigl(X(j,n),\,X(j+1,n)\bigr)$, it is also shown that $\alpha_n\rightarrow\mu$ weakly a.s. for some r.p.m. $\mu$. Now, Hill proves that, if the $\alpha_n$ are the predictives of $X$, then $X$ is not exchangeable. Therefore, in the first view, the inferrer can not have both exchangeable data and Hill's predictive distributions $\alpha_n$. On the contrary, in the second view, the inferrer can assume $X$ exchangeable, say $P(X\in\cdot)=E\bigl\{\mu^\infty(\cdot)\bigr\}$ for some r.p.m. $\mu$. Moreover, he/she can take $\beta_n$ as Hill's predictive distribution and regard such $\beta_n$ as an estimator of $\mu$.
\end{ex}

Our last two examples on M1 concern quite popular frameworks.

\begin{ex}\textbf{(M1; Bayesian nonparametrics)}
In the standard model of Bayesian nonparametrics, $X$ is assumed to be exchangeable. The prior is the probability distribution of $\mu$, where $\mu$ is the r.p.m. such that $P(X\in\cdot)=E\bigl\{\mu^\infty(\cdot)\bigr\}$. The usual procedure is to select a prior and then to calculate (or to approximate) the corresponding posterior, meant as the conditional distribution of $\mu$ given the observed data. Hence, it seems natural to adopt the second view about $X$. Nevertheless, sometimes, the inferrer also assigns a sequence $(\beta_n)$ of ``predictive distributions" and speaks of ``the prior induced by $(\beta_n)$"; see e.g. \cite{FPBRA,FPJB,HMW}. This can create some ambiguity. Indeed:

\medskip

\begin{itemize}

\item As already noted, in the second view, the $\beta_n$ are not the predictive distributions. Hence, $(\beta_n)$ does not induce any prior.

\item Suppose now that the inferrer adopts the first view and selects the predictives $\alpha_n=\beta_n$. Then, $X$ is not necessarily exchangeable. Indeed, in the first view, $P$ is the law induced by $(\alpha_n)$ via Ionescu-Tulcea theorem, and there is no reason for $P$ to be exchangeable (even if the $\alpha_n$ are invariant under permutations of $(X_1,\ldots,X_n)$). See \cite{FLR} and Section 1.2 of \cite{STS25}.

\end{itemize}

\medskip

\noindent To sum up, the phrase ``the prior induced by $(\beta_n)$" makes sense only in the first view, provided $\alpha_n=\beta_n$ and the law induced by $(\beta_n)$ is exchangeable. This is still true, incidentally, even if ``exchangeable" is replaced by ``partially exchangeable".

\end{ex}

\begin{ex}\textbf{(M1; Predictive resampling)}
This is a flexible method for making inference on a random parameter introduced by Fong, Holmes and Walker in \cite{FHW}; see also \cite{FPJB}. We just give a very brief sketch of the method, since our only goal is to place it under the first or the second view. At time $n$, after observing $(X_1,\ldots,X_n)=x$ for some $x\in\mathbb{R}^n$, the inferrer aims to learn about a random parameter $\theta=f(x,X_{n+1},\ldots)$ where $f$ is a known measurable function on $\mathbb{R}^\infty$. To this end, after fixing $N>n$, the probability distribution of $(X_{n+1},\ldots,X_N)$ is modeled through the predictives $\alpha_n^{(x)},\alpha_{n+1}^{(x)},\ldots,\alpha_{N-1}^{(x)}$. Finally, by suitably sampling these predictives, one obtains pseudo-data which are exploited to estimate $\theta$. We omit the details but raise the question: ``Under which view (the first or the second) the data generating mechanism is meant ?" Our feeling is that predictive resampling should be meant according to the first view, but it is not so clear. For instance, at page 1364 of \cite{FHW}, one reads: ``The Bayesian bootstrap is equivalent to the martingale posterior if we define our sequence of predictive probability distribution functions to be the sequence of empirical distribution
functions". In our notation, this means $\alpha_n=\mu_n$ for every $n$. However, as noted in Example \ref{emp}, $\alpha_n=\mu_n$ makes sense in the second view but not in the first. Indeed, in the first view, $\alpha_n=\mu_n$ yields $X_1=X_2=X_3=\ldots$ See also pages 1408 and 1414 of \cite{FHW}.
\end{ex}

We next turn to M2. In the subsequent examples, we assume that a sequence $(\alpha_n)$ of predictives has been selected and $P$ is the probability measure induced by $(\alpha_n)$ via Ionescu-Tulcea theorem.

We begin with two obvious (and possibly naive) examples, which still emphasize some possible shortcomings of M2.

\begin{ex}\textbf{(M2; i.i.d. limit sequence $Z$)}\label{xc67nm9}
It is not hard to build predictives $\alpha_n$ such that $\alpha_n\rightarrow\mu$ weakly a.s., where $\mu$ is a fixed (i.e., non-random) probability measure on $\mathcal{B}$; see e.g. Example 4.3 of \cite{ECP26}. In this case, the limit sequence $Z$ is i.i.d. with common distribution $\mu$. However, under the law $P$ induced by $(\alpha_n)$, the data sequence $X$ may be anything but i.i.d.
\end{ex}

\begin{ex}\textbf{(M2; c.i.d. data sequence $X$)}\label{bnk08u}
Suppose that $X$ is c.i.d. but not exchangeable. Then, $\alpha_n\rightarrow\mu$ weakly a.s., for some r.p.m. $\mu$, but the exchangeable sequence $Z$ driven by $\mu$ may be very different from $X$, even asymptotically. For instance, if $X_n$ converges in probability, one obtains $Z_n=Z_1$ a.s. for all $n$; see e.g. Example 1.2 of \cite{AOPCID}. Or else, the predictives $\alpha_n$ may fail to be invariant under permutations of $(X_1,\ldots,X_n)$. For instance, this happens when the $\alpha_n$ are the copula based predictive distributions introduced in \cite{HMW}. Another simple example where the $\alpha_n$ are not invariant is $\alpha_n=q^n\alpha_0+(1-q)\,\sum_{i=1}^nq^{n-i}\delta_{X_i}$, where $q\in (0,1)$ is a given constant. Finally, as noted in forthcoming Example \ref{bg67j}, $X$ and $Z$ may behave very differently under the CLT.
\end{ex}

The next two examples focus on the asymptotics of $X$ and $Z$. Among other things, Example \ref{cs22} shows that asymptotic exchangeability does not imply even the SLLN, while Example \ref{bg67j} deals with asymptotic credible intervals for a (random) parameter.

\begin{ex}\textbf{(M2; Asymptotic normality)}\label{cs22}
Let
$$X_n=\frac{\sum_{i=1}^nU_i}{\sqrt{n}}+V_n,$$
where $U$ and $V$ are independent i.i.d. sequences with $U_1$ and $V_1$ standard normal. Arguing as in Example 4.6 of \cite{ECP26}, it can be shown that:

\medskip

\begin{itemize}

\item $X$ is asymptotically exchangeable; in fact, $(X_{n+1},X_{n+2},\ldots)\overset{d}\longrightarrow Z$ where
$$Z=\bigl(U_1+V_1,\,U_1+V_2,\,U_1+V_3,\,\ldots\bigr);$$

\item $\alpha_n$ does not converge weakly a.s. (it does not converge either weakly in probability).

\end{itemize}

\medskip

\noindent Despite $X\sim Z$, however, $X$ and $Z$ are quite different, even asymptotically. We support this claim with three facts (but others could be mentioned). Firstly, it can be shown that $\overline{X}_n=\frac{1}{n}\,\sum_{i=1}^nX_i$ fails to converge in probability. Hence, $Z$ satisfies the SLLN while $X$ does not. Secondly, $X$ and $Z$ have singular probability distributions. Thirdly, the predictive distributions of $Z$ converge (weakly a.s.) while those of $X$ do not.

Some more evidence that $(X_{n+1},X_{n+2},\ldots)$ is far from being exchangeable, even if $n$ is large, comes from the covariances. If $(X_{n+1},X_{n+2},\ldots)$ was (approximately) exchangeable, then $\text{Cov}(X_{n+i},X_{n+j})$ would be (approximately) constant with respect to $i$ and $j$. On the contrary, a straightforward calculation shows that
$$\text{Cov}(X_{n+i},X_{n+j})=\sqrt{\frac{n+i}{n+j}}\quad\text{ for all }n\ge 1\text{ and }j>i\ge 1.$$
Therefore, whatever $n$ is, $\text{Cov}(X_{n+i},X_{n+j})$ swings a lot and can be made arbitrarily close to 0 or to 1.

\end{ex}

\medskip

\begin{ex}\textbf{(M2; Central limit theorems)}\label{bg67j}
Let $Z$ be an exchangeable sequence such that $E(Z_1^2)<\infty$. Define $\overline{Z}_n=\frac{1}{n}\,\sum_{i=1}^nZ_i$ and $\theta(Z)\overset{a.s.}=\lim_n\overline{Z}_n$. Then, $\sqrt{n}\,\bigl\{\overline{Z}_n-\theta(Z)\bigr\}$ converges stably (in particular, in distribution) to a mixture of centered Gaussian distributions; see e.g. Theorem 3.1 of \cite{AOPCID}. In addition,
\begin{gather}\label{bh78g}
n\,\frac{\overline{Z}_n-\theta(Z)}{\sqrt{\sum_{i=1}^n(Z_i-\overline{Z}_n)^2}}\,\overset{stably}\longrightarrow\,N(0,1)
\end{gather}
provided $\lim_n\frac{1}{n}\,\sum_{i=1}^n(Z_i-\overline{Z}_n)^2>0$ a.s. Condition \eqref{bh78g} allows to build asymptotic credible intervals for $\theta(Z)$. Having noted this fact, suppose $E(X_1^2)<\infty$, $X$ is asymptotically exchangeable, and the object of inference is $\theta(X)\overset{a.s.}=\lim_n\overline{X}_n$ (where the a.s. existence of the limit is an assumption). Then, condition \eqref{bh78g} generally fails if $Z$ is replaced by $X$. In particular, $\sqrt{n}\,\bigl\{\overline{X}_n-\theta(X)\bigr\}$ may fail to converge in distribution even if $X$ is c.i.d. and $E(X_1^2)<\infty$; see e.g. Example 3.2 of \cite{AOPCID}. Therefore, to make asymptotic inference on $\theta$, using $X$ or $Z$ is quite different even though $X\sim Z$.
\end{ex}

Examples \ref{xc67nm9}-\ref{bg67j} show that $X$ and $Z$ may have very different asymptotic behaviors even if $X\sim Z$. The next example deals with a different issue.

\begin{ex}\textbf{(M2; Choice between different predictives)}
Let $(\alpha_n)$ and $(\alpha_n^*)$ be two sequences of predictives such that $d(\alpha_n,\alpha_n^*)\rightarrow 0$. Here, $d$ is any suitable distance between probability measures, such as the bounded Lipschitz metric or the Prohorov distance. Suppose the inferrer has to choose between $(\alpha_n)$ and $(\alpha_n^*)$. This choice may seem minor. After all, since $d(\alpha_n,\alpha_n^*)\rightarrow 0$, either $\alpha_n$ and $\alpha_n^*$ both converge weakly to the same r.p.m. $\mu$ or both fail to converge. Moreover, in the former case, the exchangeable limit sequence $Z$ is driven by $\mu$. So, intuitively, using $(\alpha_n)$ or $(\alpha_n^*)$ should be more or less the same. But it is not necessarily so.

To see this, denote by $P^*$ the probability measure induced by $(\alpha_n^*)$ and recall that $P$ is the probability measure induced by $(\alpha_n)$. Define
$$\alpha_n=\mu_n\quad\text{and}\quad\alpha_n^*=\frac{c\,\alpha_0+n\,\mu_n}{n+c}$$
where $c>0$ is a constant. Then, $d(\alpha_n,\alpha_n^*)\rightarrow 0$ and $X$ is exchangeable under both $P$ and $P^*$. However, under $P$, one obtains $X_n=X_1$ a.s. for each $n\ge 1$; see Example \ref{emp}. On the contrary, under $P^*$, $X$ is a Dirichlet sequence. A similar situation occurs if $\alpha_n=\mu_n$ and $\alpha_n^*$ is Hill's predictive distribution, as defined in Example \ref{xf6yh8}. In this case, $d(\alpha_n,\alpha_n^*)\rightarrow 0$ under mild conditions on $X$. However, $X_n=X_1$ a.s. for all $n$ under $P$, while $X$ is asymptotically exchangeable (even if not exchangeable) under $P^*$.
\end{ex}

\medskip

We conclude the discussion of M2 with a remark. Obviously, to get around our criticism, one could make the notion of asymptotic exchangeability more demanding. We briefly mention two possibilities. Denote by $\nu_n$ the probability distribution of the shifted sequence $(X_{n+1},X_{n+2},\ldots)$. Recall also that, if $\gamma$ and $\nu$ are probability measures on any $\sigma$-field $\mathcal{E}$, their total variation distance is
$$\norm{\gamma-\nu}=\sup_{A\in\mathcal{E}}\,\abs{\gamma(A)-\nu(A)}.$$
Then, two possible strengthenings of asymptotic exchangeability (as defined in Section \ref{bhy67c}) are:

\medskip

\begin{itemize}

\item[(i)] $\norm{\alpha_n-\mu}\overset{a.s.}\longrightarrow 0$ for some r.p.m. $\mu$;

\item[(ii)] $\norm{\nu_n(\cdot)-P(Z\in\cdot)}\longrightarrow 0$ for some exchangeable sequence $Z$.

\end{itemize}

\medskip

\noindent Condition (i) is appealing and possibly deserves some investigations. A result related to (i) is in \cite{AOP13}. In turn, condition (ii) is probably too strong and rarely satisfied in real problems. However, under (ii), most of our discussion on M2 disappears. Moreover, condition (ii) admits a nice characterization. In fact, (ii) holds if and only if $P(X\in A)=P(Z\in A)$ for each event $A$ in the tail $\sigma$-field $\bigcap_n\sigma(X_{n+1},X_{n+2},\ldots)$.

\end{document}